\documentclass [leqno] {amsart}
  
\usepackage {amssymb}
\usepackage {amsthm}
\usepackage {amscd}

\usepackage[all]{xy}
\usepackage {hyperref}

\begin{document}
\bibliographystyle{alpha}
\theoremstyle{plain}
\newtheorem{proposition}[subsubsection]{Proposition}
\newtheorem{lemma}[subsubsection]{Lemma}
\newtheorem{corollary}[subsubsection]{Corollary}
\newtheorem{thm}[subsubsection]{Theorem}
\newtheorem{introthm}{Theorem}
\newtheorem*{thm*}{Theorem}
\newtheorem{conjecture}[subsubsection]{Conjecture}
\newtheorem{fails}[subsubsection]{Fails}

\theoremstyle{definition}
\newtheorem{definition}[subsubsection]{Definition}
\newtheorem{notation}[subsubsection]{Notation}
\newtheorem{condition}[subsubsection]{Condition}
\newtheorem{example}[subsubsection]{Example}
\newtheorem{claim}[subsubsection]{Claim}
\newtheorem{question}[subsubsection]{Question}

\theoremstyle{remark}
\newtheorem{remark}[subsubsection]{Remark}

\numberwithin{equation}{subsection}

\newcommand{\eq}[2]{\begin{equation}\label{#1}#2 \end{equation}}
\newcommand{\ml}[2]{\begin{multline}\label{#1}#2 \end{multline}}
\newcommand{\mlnl}[1]{\begin{multline*}#1 \end{multline*}}
\newcommand{\ga}[2]{\begin{gather}\label{#1}#2 \end{gather}}
\newcommand{\mat}[1]{\left(\begin{smallmatrix}#1\end{smallmatrix}\right)}

\newcommand{\arir}{\ar@{^{(}->}}
\newcommand{\aril}{\ar@{_{(}->}}
\newcommand{\are}{\ar@{>>}}

\newcommand{\xr}[1] {\xrightarrow{#1}}
\newcommand{\xl}[1] {\xleftarrow{#1}}
\newcommand{\lra}{\longrightarrow}
\newcommand{\inj}{\hookrightarrow}

\newcommand{\mf}[1]{\mathfrak{#1}}
\newcommand{\mc}[1]{\mathcal{#1}}

\newcommand{\CH}{{\rm CH}}
\newcommand{\Gr}{{\rm Gr}}
\newcommand{\codim}{{\rm codim}}
\newcommand{\cd}{{\rm cd}}
\newcommand{\Spec} {{\rm Spec}}
\newcommand{\supp} {{\rm supp}}
\newcommand{\Hom} {{\rm Hom}}
\newcommand{\End} {{\rm End}}
\newcommand{\id}{{\rm id}}
\newcommand{\Aut}{{\rm Aut}}
\newcommand{\sHom}{{\rm \mathcal{H}om}}
\newcommand{\Tr}{{\rm Tr}}
\newcommand{\coker}{{\rm coker}}


\renewcommand{\P} {\mathbb{P}}
\newcommand{\Z} {\mathbb{Z}}
\newcommand{\Q} {\mathbb{Q}}
\newcommand{\C} {\mathbb{C}}
\newcommand{\F} {\mathbb{F}}
\newcommand{\R}{\mathbb{R}}

\newcommand{\OO}{\mathcal{O}}
\newcommand{\Ext}{{\rm Ext}}

\title{On $p$-adic periods for mixed Tate motives over a number field}

\author{Andre Chatzistamatiou}
\address{Fachbereich Mathematik \\ Universit\"at Duisburg-Essen \\ 45117 Essen, Germany }
\email{a.chatzistamatiou@uni-due.de}
\author{S\.Inan \"Unver}
\address{Mathematics Department \\ Ko\c{c} University \\ 34450, Istanbul, Turkey}
\email{sunver@ku.edu.tr}

\thanks{This work has been supported by the SFB/TR 45 ``Periods, moduli spaces and arithmetic of algebraic varieties''}

\begin{abstract}
For a number field, we have a Tannaka category  of mixed Tate motives at our disposal. 
We construct $p$-adic points of the associated Tannaka group by using $p$-adic
Hodge theory. Extensions of two Tate objects yield functions on the Tannaka group, and we 
show that evaluation at our $p$-adic points is essentially given by the
inverse of the Bloch-Kato exponential map. 
\end{abstract}
 
\maketitle

\tableofcontents

\section*{Introduction}

For a number field $E$, one has
an abelian category
of mixed Tate motives $MT(E)$ \cite{DG}. A mixed Tate motive comes equipped with a weight
filtration $W$, and the associated graded pieces are sums of Tate objects.
There is a natural fibre functor $\omega$ defined by
\begin{equation*}\label{equation-intro-fibre-functor}
  \omega(M)=\bigoplus_{n\in \Z} \Hom(\Q(n),{\rm gr}_{-2n}^{W}(M));
\end{equation*}
we denote by $G_{\omega}$ the corresponding Tannaka group. 

If $\OO$ denotes the ring of integers of $E$ and
$x\in \Spec(\OO)$ is a closed point, then Deligne and Goncharov construct a Tannaka subcategory 
$MT(\OO_x)$ of $MT(E)$ consisting of  motives which are unramified at $x$ \cite[1.6]{DG}. We will denote
its group of tensor automorphisms by $G_x$. 

To a mixed Tate motive $M$ we can attach its $p$-adic realization $M_p$ which 
is a representation of the Galois group of $E$ with coefficients in $\Q_p$.
If the point $x$ lies over the prime $p$, then we can restrict in order to obtain
 a $p$-adic representation $M_{\iota,p}$ for the 
Galois group of the completion $E_x$ at $x$.
We will show that $M_{\iota,p}$ is always semistable. Furthermore, $M_{\iota,p}$ is 
crystalline if and only if $M$ is unramified at $x$, i.e.~$M\in MT(\OO_x)$  (Theorem \ref{thm-unramified=crystalline}). In fact, $p$-adic representations attached to mixed Tate motives
are contained in an abelian subcategory which admits a fibre functor $\tau$ 
similar to $\omega$. 
Denoting by $H_{\tau}$ the 
corresponding Tannaka group over $\Q_p$, $p$-adic realization yields a group 
homomorphism 
$$
H_{\tau}\xr{} G_{\omega}\otimes_{\Q}\Q_p. 
$$ 

The main purpose of this paper is to construct an $E_{x,st}$-valued 
point $\eta_{st}$ of $H_{\tau}$, where $\Spec(E_{x,st})$ is a $1$-dimensional
affine space over the field $E_x$. The $E_x$-valued points of $\Spec(E_{x,st})$
correspond naturally to the extensions of the canonical logarithm 
$\log:\OO^{\times}_{E_x}\xr{} E$ to $E^{\times}$. Therefore, any choice of such 
an extension induces via $\eta_{st}$ an $E_{x}$-valued point of $H_{\tau}$ and  
$G_{\omega}$. For the Tannaka subcategory of crystalline representations 
the picture is simpler: if $H_{\tau,cris}$ denotes their Tannaka group and 
$\pi:H_{\tau}\xr{} H_{\tau,cris}$ is the projection, then 
$\pi\circ \eta_{st}$ factors through  
$\Spec(E_x)$ and we obtain an $E_x$-valued point $\eta$ of $H_{\tau,cris}$. 
We denote by $\eta^{ur}_x$ the image of $\eta$ in $G_x$. 

To state our main theorem, we need to recall how extensions $M$ of $\Q(0)$
by $\Q(n)$ in $MT(\OO_x)$ give rise to functions on $G_x$ for $n\geq 1$. The 
natural isomorphisms $\alpha:\Q\xr{} \Hom(\Q(n),{\rm gr}^{W}_{-2n}M)$ and 
$\beta:\Hom(\Q(0),{\rm gr}^{W}_{0}M)\xr{} \Q$ induce elements 
$\alpha^{-1}\in \omega(M)^{\vee}$ and $\beta^{-1}\in \omega(M)$; we set 
$M(\eta^{ur}_x)=\alpha^{-1}(\eta^{ur}_x\cdot \beta^{-1}(1))$.

\begin{thm*}[Theorem~\ref{thm-main}]
For all $n\geq 1$, the map 
$$
\Ext^1_{MT(\OO_x)}(\Q(0),\Q(n)) \xr{} E_x, \quad M\mapsto M(\eta^{ur}_x), 
$$
is the composition of the p-adic realization 
$$
\Ext^1_{MT(\OO_x)}(\Q(0),\Q(n)) \xr{} \Ext^1 _{crys}(\Q_p(0),\Q_p(n))
$$ 
and 
the inverse of the Bloch-Kato exponential map (\ref{equation-Bloch-Kato-exp-1}).
\end{thm*}
 
\section{Filtered $\phi$-modules and mixed Tate filtered  $\phi$-modules}

\subsection{Mixed Tate filtered $\phi$-modules}

\subsubsection{}\label{section on mixed tate phi modules}
Let $K$ be a $p$-adic field with residue field $k$, i.e.~${\rm char}(K)=0$, 
$K$ is complete with respect to a fixed discrete valuation  and the 
residue field $k$ is perfect of characteristic $p$. 
Let $W(k)$ be the ring of Witt vectors of $k$, $\sigma:W(k) \to W(k)$ the 
Frobenius lift  and  $K_{0}$ the field of fractions of $W(k).$

\subsubsection{} \label{section:(phi,N)-modules}
We denote by $MF^{\phi} _{K}$ the category of filtered $\phi$-modules, i.e.~the 
objects are triples $(M,\phi,F)$, where $(M,\phi)$ is  an isocrystal over
$K_0$ and $F$ is a descending, exhaustive and separated filtration  
on $M_K=M\otimes_{K_0} K$. We denote by $MF^{\phi,N} _{K}$ the category
of filtered $(\phi,N)$-modules, i.e.~objects are tuples $(M,\phi,N,F)$
with $(M,\phi,F)\in MF^{\phi}_{K}$ and $N:M\xr{} M$ is a $K_0$-linear
endomorphism such that $N\phi=p\phi N$. We consider $MF^{\phi} _{K}$
as full subcategory of $MF^{\phi,N} _{K}$ via the functor 
$
(M,\phi,F)\mapsto (M,\phi,0,F).
$ 

The Dieudonn\'{e}-Manin classification \cite[II, \textsection 4.1]{man}  implies, by descent,  
that every isocrystal $(M,\phi)$ over $K_0$ admits a slope decomposition
$$
M=\bigoplus_{\lambda\in \Q} M_{\lambda},
$$
with $\phi(M_{\lambda})=M_{\lambda}$ and $(M_{\lambda},\phi\mid_{M_{\lambda}})$ is isoclynic
of slope $\lambda$. From the relation $N\phi=p\phi N,$ it follows that $N(M_{\lambda}) \subseteq M_{\lambda -1}$. In the following, we will use the notation:
$$
M_{\leq \lambda}:=\bigoplus_{\substack{\lambda'\in \Q\\ \lambda' \leq \lambda}} M_{\lambda'},\quad  
M_{\geq \lambda}:=\bigoplus_{\substack{\lambda'\in \Q\\ \lambda' \geq \lambda}} M_{\lambda'}.
$$

\begin{definition}\label{definition-MTK}
We say that an object $(M,\phi,F)\in MF^{\phi}_{K}$ is a \emph{mixed Tate filtered  $\phi$-module}
if the following properties are satisfied:
\begin{enumerate}
\item There is an isomorphism of $\phi$-modules 
$$
(M,\phi)\cong \bigoplus_{i\in I} (K_0,p^{n_i}\sigma),
$$
for some index set $I$, and $n_i\in \Z$. 
\item For all $i\in \Z$ the natural map
$$
F^iM_K \xr{} M_{\geq i}\otimes_{K_0}K
$$
is an isomorphism.
\end{enumerate}
We say that $(M,\phi,N,F)\in MF^{\phi,N}_{K}$ is a \emph{mixed Tate filtered $(\phi,N)$-module}
if $(M,\phi,F)$ is a mixed Tate filtered $\phi$-module.

We denote by $MT^{\phi}_K$ (resp.~$MT^{\phi,N}_K$) the full subcategory of $MF^{\phi}_{K}$ (resp.~$MF^{\phi,N}_{K}$)  with mixed Tate filtered  $\phi$-modules (resp.~$(\phi,N)$-modules) as objects. The categories $MT^{\phi}_K$ and $MT^{\phi,N}_{K}$ are additive. Again, we consider $MT^{\phi}_{K}$ as full subcategory of $MT^{\phi,N}_{K}$.
\end{definition}

For $(M,\phi,N,F)\in MT^{\phi,N}_K$,
it follows from Property (1) that all the slopes of $(M,\phi)$ are integers. From Property (2)
we conclude that the Hodge polygon of $(M_K,F)$ equals the Newton polygon of 
$(M,\phi)$.

\begin{definition}(Tate objects)\label{definition-Tate-objects}
Let $n\in \Z$ be an integer. We define the \emph{Tate object}  $K(n)\in MT^{\phi}_K$ 
by 
$$
K(n):=(K_0,p^{-n}\sigma,F),
$$
with $F$ defined by 
$$
F^j=\begin{cases} K & \text{if $j\leq -n$,} \\ 0 & \text{if $j>-n$.} \end{cases}
$$
\end{definition}

\begin{definition}(Weight filtration)\label{definition-weight-filtration}
  Let $(M,\phi,N,F)\in MT^{\phi,N}_K$. Let $i\in \Z$ be an integer. We define 
an object $W_{2i}(M,\phi,N,F)$ in $MF^{\phi}_K$ by 
$$
W_{2i}(M,\phi,N,F):=
(M_{\leq i}, \phi\mid_{M_{\leq i}}, N\mid_{M_{\leq i}},F\cap M_{\leq i}).
$$
We define an object ${\rm gr}^W_{2i}(M,\phi,N,F)$ in $MT^{\phi}_K$ by 
$$
{\rm gr}^W_{2i}(M,\phi,N,F):=(M_i,\phi\mid_{M_i},\tilde{F}),
$$
where $\tilde{F}$ is defined as follows:
$$
\tilde{F}^iM_i=M_i, \quad \tilde{F}^{i+1}M_i=0.
$$ 
\end{definition}
Note that $N(M_{\leq i})\subset M_{\leq i-1}$ and $N\mid_{M_{\leq i}}$ is well-defined.

\begin{proposition}\label{proposition-properties-weight-filtration}
 Let $(M,\phi,N,F)\in MT^{\phi,N}_K$ and $i\in \Z$. The following statements hold.
 \begin{enumerate}
 \item The object $W_{2i}(M,\phi,N,F)$ is contained in $MT^{\phi,N}_K$.
 \item There is an exact sequence 
   \begin{equation}\label{equation-weight-exact-sequence}
0\xr{} W_{2(i-1)}(M,\phi,N,F) \xr{} W_{2i}(M,\phi,N,F) \xr{} {\rm gr}^W_{2i}(M,\phi,N,F)
\xr{} 0.
   \end{equation}
 \end{enumerate}
 \begin{proof}
It is sufficient to prove the statement for $(M,\phi,0,F)$, i.e.~for objects in 
$MT^{\phi}_K$.

   For (1). It is obvious that 
$$
W_{2i}(W_{2(i+1)}(M,\phi,F))=W_{2i}(M,\phi,F),
$$
for all $(M,\phi,F)$. Therefore we may reduce to the case 
$$
W_{2(i+1)}(M,\phi,F)=(M,\phi,F).
$$ 
In this case 
$M=M_{\leq i}\oplus M_{i+1}$, and we have to prove that for all $j\in \Z$
the map 
$$
F^j\cap (M_{\leq i}\otimes_{K_0}K) \xr{} (M_{\leq i})_{\geq j}\otimes_{K_0}K
$$
is an isomorphism. Since $(M,\phi,F)$ is an object in $MT_K ^{\phi}$, 
the map is injective. In particular, the map is an isomorphism 
for all $j\geq i+1$. 

We need to show the surjectivity for $j\leq i$. By assumption, for 
every $m\in (M_{\leq i})_{\geq j}\otimes_{K_0}K$ there exists a preimage 
$m'\in F^jM_K$. By definition, the projection of $m'$ to $M_{i+1}\otimes_{K_0}K$ vanishes,
thus $m'\in F^j\cap (M_{\leq i}\otimes_{K_0}K)$.

For (2). There is an obvious morphism 
$W_{2(i-1)}(M,\phi,F) \xr{} W_{2i}(M,\phi,F)$ in $MT_K ^{\phi}.$  The morphism 
$W_{2i}(M,\phi,F) \xr{} {\rm gr}^W_{2i}(M,\phi,F)$ is defined by the 
projection $M_{\leq i}\xr{} M_i$. Since 
$F^{i+1}\cap (M_{\leq i}\otimes_{K_0}K)=0$, the projection is compatible with
the filtrations. Therefore the sequence (\ref{equation-weight-exact-sequence})
is well-defined. 

In order to prove that the sequence is exact we need to show that 
it is an exact sequence of $\phi$-modules and an exact sequence of 
filtered $K$-vector spaces. The first statement is obvious. For the
second statement we note that all members in the sequence 
(\ref{equation-weight-exact-sequence}) are objects in $MT^{\phi}_K$, thus
the Hodge polygons equal the Newton polygons. In particular, 
$$
\dim (F^j\cap M_{\leq i})= \dim (F^j\cap M_{\leq i-1})+\dim \tilde{F}^j,
$$ 
for all $j\in \Z$. This immediately implies the claim.    
\end{proof}
\end{proposition}

\begin{corollary}\label{corollary-Tate-is-wekly-admissible}
The category $MT^{\phi,N}_K$ is contained in the category of weakly admissible
filtered $(\phi,N)$-modules. 
\begin{proof}
We use the fact that weakly admissible filtered $(\phi,N)$-modules are stable
under extensions. Therefore the claim follows from Proposition 
\ref{proposition-properties-weight-filtration} provided 
we prove that ${\rm gr}^W_{2i}(M,\phi,N,F)$ is weakly admissible
for all $(M,\phi,N,F)\in MT^{\phi,N}_K$ and all $i\in \Z$. 
By Definition \ref{definition-weight-filtration}, ${\rm gr}^W_{2i}(M,\phi,N,F)$
is isomorphic to a direct sum of Tate objects $K(-i)$. 
Since Tate objects are (weakly) admissible, we are done.
\end{proof}
\end{corollary}

In contrast to the category $MF_{K}^{\phi,N}$, the category of weakly admissible
filtered $(\phi,N)$-modules $MF_{K} ^{\phi,N,wa}$ is an abelian category.

\begin{proposition}\label{proposition-MTK-abelian}
Let $f:(M,\phi_M,N_M,F_M)\xr{} (M',\phi_{M'},N_{M'},F_{M'})$ be a morphism in 
$MT^{\phi,N}_K$. We denote by $\ker(f)$ and 
$\coker(f)$ the kernel of $f$ and the cokernel of $f$ in $MF_{K}^{\phi,N,wa}$,
respectively. Then $\ker(f)$ and $\coker(f)$ are contained in $MT^{\phi,N}_K$.
In particular, $MT^{\phi,N}_K$ is an abelian category.
\begin{proof}
First, consider the full subcategory $\mc{C}$ of isocrystals over $K_0$
with objects $(M,\phi)$ such that there exists an isomorphism
$$
(M,\phi) \cong \bigoplus_{i\in I} (K_0,p^{n_i}\sigma).
$$   
It is easy to see that $\mc{C}$, as subcategory of the category of isocrystals,
contains all the kernels and cokernels of morphisms in $\mc{C}$.

We denote by $f_0$ the induced morphism $(M,\phi_M)\xr{} (M',\phi_{M'})$. 
Then 
$$
\ker(f)=(\ker(f_0),\phi\mid_{\ker(f_0)},N\mid_{\ker(f_0)},F\cap (\ker(f_0)\otimes_{K_0}K)).
$$
We know that $\ker(f_0)\in \mc{C}$ and thus satisfies Property (1)
of Definition \ref{definition-MTK}. It remains to show that
$$
F_M^i\cap (\ker(f_0)\otimes_{K_0}K) \xr{} \ker(f_0)_{\geq i}\otimes_{K_0}K 
$$
is an isomorphism. We have a commutative diagram 
$$
\xymatrix
{
0\ar[r]
&
F_M^i\cap (\ker(f_0)\otimes_{K_0}K) \ar[r] \ar[d]
&
F_M^i \ar[r] \ar[d]^{\cong}
&
F_{M'}^i \ar[d]^{\cong}
\\
0\ar[r]
&
\ker(f_0)_{\geq i}\otimes_{K_0}K\ar[r]
&
M_{\geq i}\otimes_{K_0}K \ar[r]
&
M'_{\geq i}\otimes_{K_0}K.
}
$$
Moreover, both rows are exact, which implies Property (2) of Definition \ref{definition-MTK}.

The claim for the cokernel follows dually.  
\end{proof}
\end{proposition}

\subsubsection{}  
The categories $MT^{\phi,N}_K$ and $MT^{\phi}_K$  are  $\Q_p$-linear rigid $\otimes$-categories.

\begin{lemma} \label{lemma-MTK-Tannaka}
The functor  
\begin{equation}\label{equation-fibre-functor}
\tilde{\omega}:MT^{\phi,N}_{K}  \xr{} \text{($\Q_p$-vector spaces)}, 
\quad  (M,\phi,N,F)\mapsto \bigoplus_{n\in \Z} \tilde{\omega}_n(M,\phi,F),  
\end{equation}
with 
$$
\tilde{\omega}_n(M,\phi,F)= \Hom_{MT^{\phi}_K}(K(n),{\rm gr}^W_{-2n}(M,\phi,F)),
$$
is a fibre functor. In particular, $(MT^{\phi,N}_{K},\tilde{\omega})$ and 
$(MT_{K}^{\phi} ,\tilde{\omega})$ are  Tannaka categories.
\begin{proof}
 It is easy to see that $\tilde{\omega}$ is a $\otimes$-functor. In order to see 
that $\tilde{\omega}$ is exact and faithful we will prove the existence of 
an isomorphism 
\begin{equation} \label{equation-omegaK_0-forget-isomorphism}
\tilde{\omega}_{K_0} \xr{\cong} (\gamma: (M,\phi,N,F)\mapsto M),   
\end{equation}
where $\tilde{\omega}_{K_0}(M,\phi,N,F)=\tilde{\omega}(M,\phi,N,F)\otimes_{\Q_p}K_0$ and 
$\gamma$ forgets about $\phi$, $N$ and $F$. Since $\gamma$
is exact and faithful, this will imply the claim. 

In order to construct (\ref{equation-omegaK_0-forget-isomorphism}), we observe
that there is a functorial isomorphism 
\begin{align}\label{align-omega_n-to-forget}
\Hom_{MT^{\phi}_K}(K(n),{\rm gr}^W_{-2n}(M,\phi,N,F))\otimes_{\Q_p}K_0 &\xr{} M_{-n}, \\
\phi\otimes a \mapsto a\cdot \phi(1). \nonumber  
\end{align}
\end{proof}
\end{lemma}

\begin{proposition}\label{proposition-criterion-for-MTK}
An object $(M,\phi,N,F)\in MF^{\phi,N,wa}_K$ belongs to $MT^{\phi,N}_K$ if and only if 
there exists an increasing exhaustive separated filtration $W$ by subobjects 
of $(M,\phi,N,F)$ in $MF^{\phi,N,wa}_K$ such that $W_i/W_{i-1}$ vanishes if $i$ is odd,
and is  a sum of Tate objects $K(-\frac{i}{2})$ if $i$ is even.  
\begin{proof}
For $(M,\phi,N,F)\in MT_K ^{\phi,N}$,  such a filtration exists by Definition \ref{definition-weight-filtration},
Proposition \ref{proposition-properties-weight-filtration}, and the fact that
${\rm gr}^W_{2i}(M,\phi,N,F)$ is a sum of Tate objects $K(-i)$.

Suppose now that $(M,\phi,N,F)\in MF^{\phi,N,wa}_K$ admits a filtration $W$ satisfying
the assumptions. It is easy to see that $(M,\phi)$ satisfies Property (1)
of Definition \ref{definition-MTK}.  

In general, if 
$$
0\xr{} M_1 \xr{} M \xr{} M_2 \xr{} 0
$$
is an exact sequence in $MF^{\phi,N,wa}_K$, and $M_1,M_2$ satisfy Property (2), then 
$M$ satisfies Property (2). By induction on $i$ we conclude that 
$W_i\in MT^{\phi,N}_K$ for all $i$.
\end{proof}
\end{proposition}

It is clear that any filtration as in Proposition \ref{proposition-criterion-for-MTK} has to coincide 
with the weight filtration, and that any morphism between two objects 
in $MT_{K} ^{\phi,N}$  has to be strict with respect to the 
weight filtrations on these objects.

\subsection{The crystalline logarithmic point}

\subsubsection{}
Recall from (\ref{equation-fibre-functor}) that we have a fibre functor $\tilde{\omega}$
equipping $MT_K ^{\phi}$ and $MT_{K} ^{\phi,N}$ with the structure of  Tannaka categories (Lemma \ref{lemma-MTK-Tannaka}).
Let $G_{\tilde{\omega}} $ and $G_{\tilde{\omega}} ^{st} $   denote 
the  pro-algebraic groups   which  represent
 tensor automorphisms of  $\tilde{\omega}$ on $MT_K^{\phi}$ and $MT_{K}^{\phi,N}$, 
respectively.  In other words, we have  $G_{\tilde{\omega}}=\underline{{\rm  Aut}}_{MT_{K}^{\phi}}^{\otimes}\tilde{\omega}$ and $G_{\tilde{\omega}} ^{st}=\underline{{\rm Aut}}_{MT_{K} ^{\phi,N}}^{\otimes}\tilde{\omega}$.
The goal of this section is to construct a non-trivial $K$-valued point $\eta$ of $G_{\tilde{\omega}}$.

\begin{definition}\label{definition-eta}
For $(M,\phi,F)\in MT^{\phi}_K$ we define 
$$
\eta(M,\phi,F): M_K \xr{} M_K
$$  
to be the unique endomorphism rendering the following diagram commutative:
\begin{equation}\label{diagram-eta}
\xymatrix
{
M_K=\bigoplus_{i\in \Z} M_i\otimes_{K_0}K \ar[rrr]^{\bigoplus_i \iota_i\otimes_{K_0}K} \ar[d]_{\eta(M,\phi,F)}
&
&
&
\bigoplus_{i\in \Z} M_{\geq i}\otimes_{K_0}K \ar[d]^{\left(\bigoplus_{i\in \Z} \pi_i\right)^{-1}}
\\
M_K
&
&
&
\bigoplus_{i\in \Z} F^iM_K \ar[lll]^{\sum_{i\in \Z}},
}  
\end{equation}
where $\iota_i:M_i\xr{} M_{\geq i}$ is the obvious inclusion, $\pi_i:F^{i}M_K\xr{} M_{\geq i}\otimes_{K_0}K$ is the 
projection and therefore by definition an isomorphism (Definition \ref{definition-MTK}(2)), and $\sum_{i\in \Z}$
is the sum over the obvious inclusions. 
\end{definition}

\begin{lemma}\label{lemma-eta-automorphism}
  The morphisms $\eta$ from Definition \ref{definition-eta} define a tensor automorphism of the fibre functor $\tilde{\omega}_K=\tilde{\omega}\otimes_{\Q_p}K$.
  \begin{proof}
Via the $\otimes$-isomorphism (\ref{equation-omegaK_0-forget-isomorphism}) we 
may identify $\tilde{\omega} \otimes_{\Q_p} K_0$ with the forgetful functor 
$(M,\phi,F)\mapsto M$. After tensoring with $K$ we obtain 
$\tilde{\omega}_K(M,\phi,F)=M_K$.
 
First, let us prove that $\eta(M,\phi,F)$ is an automorphism. We denote by 
$$
\eta(M,\phi,F)[i,j]:M_j\otimes_{K_0}K \xr{} M_i\otimes_{K_0}K
$$
the the composition with the inclusion $M_j\otimes_{K_0}K\xr{} M_K$ and the projection $M_{K}\xr{} M_i\otimes_{K_0}K$.
It is easy to see from the definitions that 
\begin{equation}\label{equation-eta-traingular}
\eta(M,\phi,F)[i,j]=
\begin{cases} 
0 &  \text{if $i>j$,} \\
id_{M_i} & \text{if $i=j$.}
\end{cases}  
\end{equation}
Therefore $\eta(M,\phi,F)$ is an automorphism.

Since the diagram (\ref{diagram-eta}) is functorial, $\eta$ defines a natural transformation. The compatibility 
with the tensor product is obvious.
\end{proof}
\end{lemma}

\subsubsection{} \label{section-comparison-with-Deligne-construction}
Let us explain the construction of $\eta$ in the formalism 
of \cite{D}. For $(M,\phi,F)\in MT^{\phi}_K$ there are three filtrations on $M_K$:
\begin{enumerate}
\item the weight filtration:
$$
W_iM_K=
\begin{cases} 
M_{\leq \frac{i}{2}}\otimes_{K_0} K & \text{if $i$ is even} \\
M_{\leq \frac{i-1}{2}}\otimes_{K_0} K & \text{if $i$ is odd.}
\end{cases}
$$
\item The Hodge filtration $F$.
\item The filtration 
$$
\bar{F}^iM_K:= M_{\geq i}\otimes_{K_0}K \qquad \text{for all $i\in \Z$.}
$$
\end{enumerate}
The three filtration $W,F,\bar{F}$ satisfy the condition 
$$
{\rm Gr}^p_F{\rm Gr}^q_{\bar{F}}{\rm Gr}^W_nM_K =0 \qquad \text{for $n\neq p+q,$}
$$
of \cite[\textsection1.1]{D}.
Induced by $F,\bar{F}$, we obtain maps 
\begin{align*}
  a_{F}:M_K &= \bigoplus_{i\in \Z} F^{i}\cap W_{2i} \xr{} \bigoplus_{i\in \Z} {\rm Gr}^W_{2i}M_K, \\
 a_{\bar{F}}: M_K &= \bigoplus_{i\in \Z} \bar{F}^{i}\cap W_{2i} \xr{} \bigoplus_{i\in \Z} {\rm Gr}^W_{2i}M_K,
\end{align*}
where $F^{i}\cap W_{2i}\xr{} {\rm Gr}^W_{2i}M_K$ is the natural map (and similarly 
for $\bar{F}$). We obtain a unipotent automorphism $d=a_{\bar{F}}a_F^{-1}$ of 
$\bigoplus_{i\in \Z} {\rm Gr}^W_{2i}M_K$. 

It is easy to see that we have the equality 
$$
\eta(M,\phi,F)=a_{\bar{F}}^{-1}\circ d\circ a_{\bar{F}}.
$$

\subsubsection{}
Let us see in explicit terms how $\eta$ compares the crystalline structure 
with the Hodge filtration.
For $(M,\phi,F)\in MT_K ^{\phi}$  we say that $v_1,\dots,v_d\in \tilde{\omega}(M,\phi,F)\otimes_{\Q_p} K$ is
a homogeneous basis if it is a basis of $\tilde{\omega}(M,\phi,F)\otimes_{\Q_p} K$, 
and for every $v_i$ there is an integer $n_i$ with 
$v_i\in \tilde{\omega}_{n_i}(M,\phi,F)\otimes_{\Q_p} K$; we set $\deg(v_i)=n_i$.

Recall that for all integers $i$ we have isomorphisms 
\begin{align*}
  a_i&:M_{-i}\otimes_{K_0} K \xr{\cong} \tilde{\omega}_{i}(M,\phi,F)\otimes_{\Q_p} K,\\
  b_i&:F^{-i} \cap W_{-2i}M_K \xr{\cong} \tilde{\omega}_{i}(M,\phi,F)\otimes_{\Q_p} K.
\end{align*}
The first map is the inverse of (\ref{align-omega_n-to-forget}) and the second
map is given by the composition
$$
b_i: F^{-i} \cap W_{-2i}M_K \xr{\subset} M_K \xr{\text{projection}} 
M_{-i}\otimes_{K_0}K \xr{a_i}  \tilde{\omega}_{i}(M,\phi,F)\otimes_{\Q_p} K. 
$$
For a homogeneous basis $\{v_j\}$ we set 
$$
v_j^{{\rm crys}}:=a_{\deg(v_j)}^{-1}(v_j), \qquad v_j^{{\rm Hodge}}:=b_{\deg(v_j)}^{-1}(v_j).
$$
We denote by $\{v_j^{{\rm crys},\vee}\}$ the basis dual to $\{v_j^{{\rm crys}}\}$.
By definition of $\eta$ we have 
$$
v_i^{\vee}(\eta(v_j))=v_i^{{\rm crys},\vee}(v_j^{{\rm Hodge}}).
$$

\subsubsection{} 
By Lemma \ref{lemma-eta-automorphism}, we obtain a $K$-valued point $\eta\in G_{\tilde{\omega}}(K)$; we call this point the 
\emph{logarithmic point}. Let us check that $\eta$ is not the identity. 

\begin{proposition}\label{proposition-MTK-extensions}
  Let $n\in \Z$ be an integer. 
We have 
\begin{equation}\label{equation-MTK-extensions}
{\rm Ext}^1_{MT^{\phi}_K}(K(0),K(n))\cong\begin{cases} K &\text{if $n>0$} \\ 0 &\text{if $n\leq 0$.} \end{cases}  
\end{equation}
Let 
\begin{equation}\label{equation-extension}
0\xr{} K(n) \xr{\iota} (E,\phi,F) \xr{\pi} K(0) \xr{} 0  
\end{equation}
be an extension. For $n\neq 0$, there are unique sections 
$f:E\xr{} K_0$ and $v:K_0\xr{} E$  of the underlying maps of $K_{0}$-isocrystals of $\iota$ and  $\pi$, respectively.
The isomorphism (\ref{equation-MTK-extensions}), for $n\neq 0$, is given by the formula
$$
E\mapsto f(\eta(E,\phi,F)(v(1))).
$$
\begin{proof}
First, we consider the case $n=0$.  Let $(E,\phi,F)$ be as in (\ref{equation-extension}). 
We have $F^0(E_K)=E_K$ and $F^1(E_K)=0$ by Definition \ref{definition-MTK}(2). In view 
of  Definition \ref{definition-MTK}(2) there is an isomorphism 
$(E,\phi)\cong (K_0,\sigma)\oplus (K_0,\sigma)$, thus there is a section of $\pi$ in $MT^{\phi}_K$.

For $n\neq 0$: From the slope decomposition we obtain natural sections $f,v$ as $\phi$-modules.   
If $n<0$ then $F^1E_K=\iota(K)$ which means $(E,\phi,F)=K(0)\oplus K(n)$.

For $n>0$, we can uniquely write $F^{n+1}E_K=K\langle a\cdot \iota(1) + v(1)\rangle$ with $a\in K$. Obviously, 
$$
f(\eta(E,\phi,F)(v(1)))=a.
$$ 
It is clear that $F^{n+1}E_K$ is the only invariant for extensions.
\end{proof}
\end{proposition}

\subsubsection{}
Recall that we have a fibre functor $\tilde{\omega}$ (\ref{equation-fibre-functor})
to the category of $\Q_p$-vector spaces. In the obvious way $\tilde{\omega}$
factors through the category of graded $\Q_p$-vector spaces. 
Furthermore,
we have an automorphism $\eta$ of $\tilde{\omega}_K$ (Lemma~\ref{lemma-eta-automorphism}).

\begin{definition}\label{definition-Ceta}
We define $\mc{C}_{\eta}$ to be the category of pairs $(V,\eta)$, where 
$V$ is a finite dimensional graded $\Q_p$-vector space and 
$\eta:V\otimes K\xr{} V\otimes K$ is a $K$-linear map such that for all $n\in \Z$:
\begin{equation}\label{equation-grading-and-eta}
(\eta-id)(V_n\otimes K) \subset \bigoplus_{i>n} V_i\otimes K.  
\end{equation}
Morphisms $(V_1,\eta_1)\xr{} (V_2,\eta_2)$ are $\Q_p$-linear morphisms 
$\tau:V_1\xr{} V_2$ which respect the grading 
and commute with the endomorphisms $\eta_i$, i.e.~$\eta_2\circ (\tau\otimes id_K)=(\tau\otimes id_K)\circ \eta_1$.
\end{definition}

The category $\mc{C}_{\eta}$ is a $\otimes$-category with 
$$
(V_1,\eta_1)\otimes (V_2,\eta_2)=(V_1\otimes V_2,\eta_1\otimes \eta_2).
$$

\begin{proposition} \label{proposition-MT_K-Ceta-equivalence}
  The functor 
  \begin{align*}
    \Psi:MT_K ^{\phi}&\xr{} \mc{C}_{\eta} \\
    (M,\phi,F)&\mapsto \left(\bigoplus_{n\in \Z} \tilde{\omega}_n(M,\phi,F), \eta(M,\phi,F) \right)
  \end{align*}
is an equivalence of $\otimes$-categories.
\begin{proof}
By Lemma \ref{lemma-eta-automorphism}, $\eta$ is functorial and $\Psi$ is a $\otimes$-functor. It follows from 
(\ref{equation-eta-traingular}) that 
$$(\eta-id)(\tilde{\omega}_n\otimes K)\subset \bigoplus_{i>n} \tilde{\omega}_i\otimes K.$$

We define a functor 
  \begin{align*}
    \Phi:\mc{C}_{\eta}&\xr{} MT_K  ^{\phi} \\
    (\oplus_{n\in \Z} V_n,\eta)&\mapsto 
\left(\oplus_{n\in \Z} (V_{-n}\otimes_{\Q_p} K_0,p^{-n}\otimes \sigma), F \right),
  \end{align*}
with the following filtration:
$$
F^i:=\eta\left( \bigoplus_{j\geq i} V_{-j}\otimes_{\Q_p} K  \right),
$$
for all $i$. Property (\ref{equation-grading-and-eta}) implies that $\Phi$
is well-defined. From Definition \ref{definition-eta} it easily follows
that $\Psi\circ \Phi=id_{\mc{C}_{\eta}}$.

On the other hand, we have $\Phi\circ \Psi \xr{\cong} id_{MT_K ^{\phi}}$ via
\begin{align*}
  \Phi\circ \Psi(M,\phi,F) &\xr{} (M,\phi,F) \\
  \bigoplus_{n\in \Z} \tilde{\omega}_{-n}(M,\phi,F)\otimes_{\Q_p} K_0 &\xr{\text{(\ref{align-omega_n-to-forget})}} M.
\end{align*}
\end{proof}
\end{proposition}

\begin{remark}
  Via the dictionary of Section \ref{section-comparison-with-Deligne-construction},
Proposition \ref{proposition-MT_K-Ceta-equivalence} is a variant of \cite[Proposition~1.2]{D}.
\end{remark}

\subsection{The semistable logarithmic point} 

\subsubsection{}
Let $K$ be as in \textsection\ref{section on mixed tate phi modules} with 
residue field $k$. We denote by  $\nu_K$ the valuation of $K$. 

\subsubsection{}
Recall that we have a homomorphism 
$$
[.]: k^{\times}\xr{} \OO_K^{\times}, \quad x\mapsto [x],
$$ 
by taking the Teichm\"uller lift. Denoting by 
$U_K:=\{x\in \OO_K^{\times}; x\in 1+m_K\}$ the $1$-units, we obtain a decomposition
$$
\OO^{\times}_K=k^{\times} \times U_K.
$$
The logarithm 
\begin{equation}\label{equation:logarithm-units}
\log:\OO_K^{\times} \xr{} \OO_K  
\end{equation}
is by definition trivial on the factor $k^{\times}$ and is given by 
$$
\log(u)=\sum_{n\geq 1}(-1)^{n+1} \frac{(u-1)^n}{n} \qquad \text{for all $u\in U_K$.} 
$$

\subsubsection{}
We consider $\OO^{\times} _{K,\Q}:= \OO^{\times}_K\otimes_{\Z}\Q$ and $K^{\times} _{\mathbb{Q}}:=K^{\times}\otimes_{\Z}\Q$  as $\Q$-vector spaces,
therefore we may form the symmetric algebras ${\rm Sym}_{\Q}(\OO^{\times}_{K,\Q})$
and ${\rm Sym}_{\Q}(K^{\times}_{\Q})$. The exact sequence 
$$
0\xr{} \OO_{K,\Q}^{\times} \xr{} K^{\times}_{\Q} \xr{\nu_K} \Q \xr{} 0
$$ 
implies that $\Spec({\rm Sym}_{\Q}(K^{\times} _{ \Q}))$ is a $1$-dimensional 
affine space over the scheme $\Spec({\rm Sym}_{\Q}(\OO^{\times}_{K,\Q}))$. In other words, 
for $x\in K^{\times}$ with $\nu_K(x)\neq 0$, the map 
$$
{\rm Sym}_{\Q}(\OO^{\times}_{K, \Q})[X]\xr{} {\rm Sym}_{\Q}(K^{\times}_ {\Q}), \quad X\mapsto x,
$$ 
is an isomorphism. 

The logarithm (\ref{equation:logarithm-units}) induces a ring homomorphism 
\begin{equation}\label{equation:logarithm-Sym-units}
{\rm Sym}_{\Q}(\OO^{\times}_{K,\Q}) \xr{} K,  
\end{equation}
because $K$ is torsion free.

\begin{definition}\label{definition-Kst}
  We define the $K$-algebra $K_{st}$ by 
$$
K_{st}:={\rm Sym}_{\Q}(K^{\times}_ {\Q})\otimes_{{\rm Sym}_{\Q}(\OO^{\times}_{K,\Q})} K.
$$
\end{definition}

By base change, $\Spec(K_{st})$ is a $1$-dimensional affine space over $K$. 
We have a natural logarithm 
\begin{equation}\label{equation-log-st}
\log_{st}: K^{\times} \xr{} K_{st}, \quad x\mapsto x\otimes 1.  
\end{equation}
The $K$-valued points of $\Spec(K_{st})$ admit the following description:
\begin{align}\label{equation-K-valued-points-K_st}
\Spec(K_{st})(K)&=\{\text{extensions $\log:K^{\times}\xr{} K$ of (\ref{equation:logarithm-units})}\}\\
  f&\mapsto f^*\circ \log_{st}.\nonumber
\end{align}
By an extension $\log:K^{\times}\xr{} K$ we mean a homomorphism such that the restriction
to $\OO_K^{\times}$ equals (\ref{equation:logarithm-units}). 

\subsubsection{}
The $p$-adic Hodge theory for $K$ (and fixed valuation $\nu_K$) 
depends for semistable representations on the choice of a logarithm 
$$
\log:K^{\times}\xr{} K.
$$ 
It will be important for us that our constructions do not depend on a 
particular choice, and for this we have to recall the basic constructions
of $p$-adic Hodge theory.

We denote by $R$ the ring
$$
R:=\varprojlim \OO_{\bar{K}}/p\OO_{\bar{K}},
$$ 
where the maps are given by rising to the $p$-th power $x\mapsto x^p$. 
Denoting by $C_K=\widehat{\overline{K}}$ the $p$-adic completion of 
$\overline{K}$ we have a multiplicative bijection 
$$
\varprojlim \OO_{C_K} \xr{} R,
$$
where the projective system is defined by rising to the $p$-th power again. In
other words, we can represent every element $x$ in $R$ by 
$(x^{(0)},x^{(1)},\dots)$ with $x^{(n)}\in \OO_{C_K}$ and $x^{(n-1)}=(x^{(n)})^p$.

Let $\nu_{\bar{K}}$ (resp.~$\nu_{C_K}$) be the extension of $\nu_{K}$ (resp.~$\nu_{\bar{K}}$) to $\bar{K}$ (resp.~$C_K$).
The map
$$
\nu_R: R\backslash \{0\} \xr{} \Q, \quad x\mapsto \nu_{C_K}(x^{(0)})  
$$
can be extended  to a valuation
$$ 
\nu_R:{\rm Frac}(R)^{\times} \xr{} \Q
$$
with valuation ring $R$.

Let $B_{cris}$ be the crystalline period ring; we define
$$
B_{st}=
{\rm Sym}_{\Q}({\rm Frac}(R)^{\times}\otimes_{\Z} \Q)\otimes_{{\rm Sym}_{\Q}(R^{\times}\otimes_{\Z} \Q)} B_{cris},
$$
where ${\rm Sym}_{\Q}(R^{\times}\otimes_{\Z}\Q)\xr{} B_{cris}$ is induced by the crystalline logarithm 
$$
\log_{cris}: R^{\times}\xr{} B_{cris}.
$$
Again, $R^{\times}=\bar{k}^{\times}\times (1+m_R)$; $\log_{cris}$ is trivial on $\bar{k}^{\times}$
and given by  
$$
\log_{cris}(u)=\sum_{n\geq 1}(-1)^{n+1} \frac{([u]-1)^n}{n}
$$
for $u\in 1+m_R$, where $[u]$ denotes the Teichm\"uller lift of $u$ in the Witt ring  $W(R)$ of $R.$ 
 
By construction we have a natural logarithm 
$$
\log_{st}: {\rm Frac}(R)^{\times} \xr{} B_{st}, \quad x\mapsto x\otimes 1.
$$

The ring $B_{st}$ has the following properties.
\begin{enumerate}
\item We have a ${\rm Gal}(\bar{K}/K)$-action on $B_{st}$ extending the action on $B_{cris}$.
\item We have a Frobenius map $\phi:B_{st}\xr{} B_{st}$ extending the Frobenius map on $B_{cris}$. Moreover, $$\phi\circ \log_{st}=p\log_{st}.$$
\item We have a $B_{cris}$-linear derivation $N:B_{st}\xr{} B_{st}$ such that  
$$
N(\log_{st}(x))=\nu_R(x) \qquad \text{for all $x\in {\rm Frac}(R)^{\times}$.}   
$$ 
\end{enumerate}

After choosing a logarithm 
$$
\log:K^{\times}\xr{} K,
$$
which extends (\ref{equation:logarithm-units}), we obtain a morphism of $B_{cris}$-algebras
$$
\gamma_{\log}:B_{st}\xr{} B_{dR}.
$$
The morphism depends on the choice of $\log$, and the filtration on $B_{st}$
induced by the filtration on $B_{dR}$ via $\gamma_{\log}$ depends on $\log$. 
In order to simplify the comparison between different logarithms we will restrict ourselves to logarithms
$\log$ such that $\log(K_0^{\times})\subset K_0$. In other words, we will only consider
$K_0$-valued points of $\Spec(K_{0,st})$.

\begin{proposition}\label{proposition-delta-log-log'}
  For $\log,\log'\in \Spec(K_{0,st})(K_0)$ there is a unique ring homomorphism 
$$
\delta_{\log,\log'}:B_{st}\xr{} B_{st}
$$
such that $\gamma_{\log'}\circ \delta_{\log,\log'}=\gamma_{\log}$. 
The map $\delta_{\log,\log'}$ is given by 
\begin{equation}\label{equation-delta-log-log'}
\delta_{\log,\log'}=\exp\left(\frac{\log(x)-\log'(x)}{\nu_K(x)} N\right)  
\end{equation}
for every $x\in K_0\backslash \OO_{K_0}^{\times}$.
\begin{proof}
Uniqueness follows from the fact that $\gamma_{\log}$ is injective.

Choose $\tilde{p}\in R$ with $\tilde{p}^{(0)}=p$. By definition we have
$$
\gamma_{\log}(\log_{st}(\tilde{p}))=\log_{dR}([\tilde{p}]/p)+\log(p),
$$
where $\log_{dR}$ is defined by the usual series since $[\tilde{p}]/p$ is a 
$1$-unit in $B_{dR}$. Since $\Spec(B_{st})$ is a $1$-dimensional affine 
space over $\Spec(B_{cris})$, there exists a unique morphism of $B_{cris}$-algebras 
$\delta_{\log,\log'}$ such that 
$$
\delta_{\log,\log'}(\log_{st}(\tilde{p}))=\log_{st}(\tilde{p})+\log(p)-\log'(p).
$$
Obviously, $\delta_{\log,\log'}$ satisfies $\gamma_{\log'}\circ \delta_{\log,\log'}=\gamma_{\log}$ and the equality (\ref{equation-delta-log-log'}).
\end{proof}
\end{proposition}

By using $\gamma_{\log}$ we obtain a filtration on $B_{st}$. The $p$-adic Hodge 
theory \cite[Thm.~A]{CF} asserts that the functor 
\begin{align}
D_{st,\log}:&(\text{semistable $\Q_p$-representations of ${\rm Gal}(\bar{K}/K)$})\xr{} MF^{\phi,N,w.a.}_K \\
&V\mapsto \left(B_{st}\otimes_{\Q_p}V\right)^{{\rm Gal}(\bar{K}/K)}  \nonumber
\end{align}
is an equivalence of categories. We will use the subscript $\log$ in $D_{st,\log}$
to emphasize 
the dependence on $\log$. 

Denoting by ${\rm forget}_{F}$ the functor 
${\rm forget}_{F}(M,\phi,N,F)=(M,\phi,N)$, we get 
$$
{\rm forget}_{F}\circ D_{st,\log}={\rm forget}_{F}\circ D_{st,\log'},
$$
because only the filtration depends on the embedding to $B_{dR}$. Proposition \ref{proposition-delta-log-log'}
implies that for the filtrations we have the following comparison:
\begin{equation}\label{equation-comparson-filtrations}
  F^i_{D_{st,\log'}(V)} = \exp\left(\frac{\log(x)-\log'(x)}{\nu_K(x)} N\right) F^i_{D_{st,\log}(V)},
\end{equation}
for all $i\in \Z$ and all $x\in K_0\backslash \OO^{\times}_{K_0}$.

\begin{definition}\label{definition-mixed-Tate-representation}
Let $K$ be a $p$-adic field.
We denote by $MT_{G_K}$ the full subcategory of $p$-adic representations $V$ 
of ${\rm Gal}(\bar{K}/K)$ which admit an increasing exhaustive separated filtration $W$ by subrepresentations 
of $V$ such that $W_i/W_{i-1}$ vanishes if $i$ is odd,
and is  a sum of Tate objects $\Q_p(-\frac{i}{2})$ if $i$ is even.
We call an object of   $MT_{G_K}$ a \emph{mixed Tate representation} of ${\rm Gal}(\bar{K}/K)$. 
\end{definition}

\begin{proposition} \label{proposition-mixed-Tate-rep-well-defined}
  Let $\log\in \Spec(K_{0,st})(K_0)$. Then 
$$
MT_{G_K}=D_{st,\log}^{-1}(MT^{\phi,N}_K).
$$
In particular, every mixed Tate representation is semistable.
\begin{proof}
From Proposition \ref{proposition-criterion-for-MTK} it follows that every object
in $D_{st,\log}^{-1}(MT^{\phi,N}_K)$ admits a filtration $W$ satisfying the properties
of Definition \ref{definition-mixed-Tate-representation}.

Now, suppose that $V$ is a $p$-adic representation of ${\rm Gal}(\bar{K}/K)$ which admits a filtration $W$ as in Definition \ref{definition-mixed-Tate-representation}. If we know that $V$ is semistable then clearly 
$D_{st,\log}(V)\in MT^{\phi,N}_K$ by Proposition \ref{proposition-criterion-for-MTK}, again. Therefore it suffices to prove that $V$ is semistable.

We use induction on the length of the filtration $W$ of $V$. 
If the filtration $W$ has length $\leq 1,$ semistability of $V$ follows from those of $\mathbb{Q}_{p}(n).$  
In general, let $n$ be the smallest integer such that $W_{2n}V = V.$ Then we have an exact sequence 
$$
0 \to (W_{2n-2}V)\otimes \mathbb{Q}_{p}(n) \to V \otimes \mathbb{Q}_{p}(n)  \to (V/W_{2n-2}V) \otimes \mathbb{Q}_{p}(n) \to 0.
$$
By the induction hypothesis the terms on the left and right are semistable. Moreover, since the weights of the term on the left are $\leq -2$ and the term on  the right has weight 0, we have 
\begin{align*}
F^{0} D_{dR}((W_{2n-2}V)\otimes \mathbb{Q}_{p}(n))&=0 \\  
F^{0} (D_{dR} ( (V/W_{2n-2}V)\otimes \mathbb{Q}_{p}(n) ) ) &= D_{dR} ( (V/W_{2n-2}V) \otimes \mathbb{Q}_{p}(n)).
\end{align*}
Therefore \cite[Proposition 1.28]{nek} shows that the middle term is also semistable. 
\end{proof}
\end{proposition}

Obviously, 
\begin{equation}\label{definition-tau}
\tau=\tilde{\omega}\circ D_{st,\log}  
\end{equation}
is independent of $\log$, and $(MT_{G_K},\tau)$ is a Tannaka category (by Lemma \ref{lemma-MTK-Tannaka}).

\subsubsection{}
Recall from Lemma \ref{lemma-MTK-Tannaka} that $(MT^{\phi,N}_{K},\tilde{\omega})$
is a Tannaka category. We will use the ring $K_{st}$ (Definition \ref{definition-Kst}) and $\log_{st}$ (\ref{equation-log-st}).

\begin{definition}\label{definition-eta-st}
For a logarithm $\log\in \Spec(K_{0,st})(K_0)$ and $(M,\phi,N,F)\in MT^{\phi,N}_K$
we define $\eta_{st,\log}(M,\phi,N,F)\in \End_{K_{st}}(M \otimes_{K_0}K_{st})$ by
$$
\eta_{st,\log}(M,\phi,N,F):=\exp\left(\frac{\log(x)-\log_{st}(x)}{\nu_K(x)} N\right)\eta(M,\phi,F),
$$ 
for $x\in K_0^{\times}\backslash \OO_{K_0}^{\times}$. For the definition of $\eta(M,\phi,F)$
we refer to Definition \ref{definition-eta}.
\end{definition}
Obviously, $\eta_{st,\log}$ does not depend on the choice $x\in K_0^{\times}\backslash \OO_{K_0}^{\times}$, but it depends on $\log$.

\begin{lemma}\label{lemma-eta-st}
Let $\log\in \Spec(K_{0,st})(K_0)$. The morphisms $\eta_{st,\log}$ from Definition \ref{definition-eta-st} define a tensor automorphism of the fibre functor $\tilde{\omega}_{K_{st}}=\tilde{\omega}\otimes_{\Q_p}K_{st}$. In other words, $\eta_{st,\log}\in G_{\tilde{\omega}}^{st}(K_{st})$ with $G_{\tilde{\omega}}^{st}=\underline{{\rm Aut}}_{MT_{K} ^{\phi,N}}^{\otimes}\tilde{\omega}$.
  \begin{proof}
Via the $\otimes$-isomorphism (\ref{equation-omegaK_0-forget-isomorphism}) we 
may identify $\tilde{\omega} \otimes_{\Q_p} K_0$ with the forgetful functor 
$(M,\phi,N,F)\mapsto M$. After tensoring with $K_{st}$ we obtain 
$\tilde{\omega}_{K_{st}}(M,\phi,N,F)=M\otimes_{K_0}K_{st}$. Lemma \ref{lemma-eta-automorphism}
implies that $\eta(M,\phi,F)$ is a tensor automorphism, thus it suffices to
prove the statement for 
$\exp\left(\frac{\log(x)-\log_{st}(x)}{\nu_K(x)} N\right)$. The functoriality
follows immediately. The compatibility with the $\otimes$-structure follows
from 
$$
N_{M_1\otimes M_2}=N_{M_1}\otimes 1 + 1 \otimes N_{M_2}.
$$
\end{proof}
\end{lemma}

\begin{lemma}\label{lemma-eta-st-well-defined}
The $K_{st}$-valued point 
$$
\eta_{st}=\eta_{st,\log}\circ D_{st,\log}
$$
of $\underline{{\rm Aut}}_{MT_{G_K}}^{\otimes}\tau$ is independent of the choice 
of $\log\in \Spec(K_{0,st})(K_0)$.
\begin{proof}
Let $\log,\log'\in \Spec(K_{0,st})(K_0)$ and $V\in MT_{G_K}$.
In view of (\ref{equation-comparson-filtrations}) we get 
\begin{equation}\label{equation-eta-compare-log-log'}
\eta({\rm forget}_ND_{st,\log'}(V))=\exp\left(\frac{\log(x)-\log'(x)}{\nu_K(x)} N\right)\eta({\rm forget}_ND_{st,\log}(V)),  
\end{equation}
for very $x\in K_0\backslash \OO_{K_0}^{\times}$, and ${\rm forget}_N(M,\phi,N,F)=(M,\phi,F)$. Thus
\begin{equation*}
\begin{split}
  \eta_{st,\log'}D_{st,\log'}(V)&= \exp\left(\frac{\log'(x)-\log_{st}(x)}{\nu_K(x)} N\right) \eta({\rm forget}_ND_{st,\log'}(V))\\
&= \exp\left(\frac{\log(x)-\log_{st}(x)}{\nu_K(x)} N\right) \eta({\rm forget}_ND_{st,\log}(V)) \quad \text{~by (\ref{equation-eta-compare-log-log'}),}\\
&=\eta_{st,\log}D_{st,\log}(V).
\end{split}
\end{equation*}
\end{proof}
\end{lemma}

\begin{example}
By Kummer theory any  $q \in K^{\times}$ defines an extension $V$ of the ${\rm Gal}(\bar{K}/K)$-representation $\mathbb{Q}_{p}(0)$ by  $\mathbb{Q}_{p}(1). $ This in turn gives via $D_{st,\log}$ an extension of $K(0)$ by $K(1)$ in $MT_{K}^{\phi, N}$: 
$$
0 \to K(1) \to M \to K(0)\to 0,
$$
which may be described as follows.   
The underlying $K_{0}$-space of $M$ has a basis $e_{0},e_{1}$ such that the 
following conditions are satisfied:
\begin{enumerate}
\item the action of $\phi$ is given by 
$\phi(e_{i})=p^{-i}e_{i}$ for $i=0,1,$ 
\item $e_1$ is the image of $1\in K(1)$,
\item $e_0$ maps to $1\in K(0)$.
\end{enumerate}

The filtration is given by  $F^{-1} M_{K}=M_{K},$ $F^{0} M_{K}=K\cdot \langle \log(q) e_{1}+e_{0}\rangle$ and $F^{1}M_{K}=0.$  Finally $N$ is given by $Ne_{0}=-\nu_K(q)\cdot e_{1}$ and $Ne_{1}=0.$ Then we easily compute 
$$
\eta_{st}(V)=\left(\begin{matrix} 1 & 0 \\ 
\frac{-\nu_K(q)(\log(x)-\log_{st}(x))}{\nu_K(x)} & 1 \end{matrix}\right)\cdot
\left(\begin{matrix}1 & 0 \\ \log(q) & 1 \end{matrix} \right)
$$
for the obvious basis of 
$$
\tau(V)=\Hom(\Q(0),\Q(0))\oplus \Hom(\Q(1),\Q(1)),
$$
and every $x\in K_0\backslash \OO^{\times}_{K_0}$, or, equivalently, for every $x\in K\backslash \OO^{\times}_{K}$. If $\nu_K(q)\neq 0$ then we can take $q=x$ in order to see that
$$
\eta_{st}(V)=\left(\begin{matrix} 1 & 0 \\ 
\log_{st}(q) & 1 \end{matrix}\right)
$$
holds for all $q\in K^{\times}$.
\end{example}

\subsection{Tannaka group scheme of mixed Tate filtered $\phi$-modules}

\begin{definition}\label{definition-L}
We define $\mc{L}$ to be the graded $\Q_p$-Lie algebra freely generated 
by $K^{\vee}$ (i.e.~the dual of $K$, where $K$ is considered as $\Q_p$-vector 
space)  in each degree $i>0$. In other words, $\mc{L}$ is defined by 
\begin{equation} \label{equation-definition-L}
  \Hom_{\text{(graded $\Q_p$-Lie algebras)}}(\mc{L},\mc{T})=\Hom_{\text{(graded $\Q_p$-vector spaces)}}\left(\oplus_{i\in \Z_{>0}} K^{\vee},\mc{T}\right), 
\end{equation}
for every graded $\Q_p$-Lie algebra $\mc{T}$.
Via (\ref{equation-definition-L}) we get $\Q_p$-linear maps 
$$a_i:K^{\vee}\xr{} \mc{L}_i, \quad \text{for all $i>0$.}$$
\end{definition}

Obviously, $\mc{L}$ is concentrated in positive degrees and each $\mc{L}_i$  
is a finite dimensional $\Q_p$-vector space. 

\begin{definition}
  For all $n>0$, we define a $\Q_p$-Lie algebra $\mc{L}_{\leq n}$ by 
$$
\mc{L}_{\leq n} = \mc{L}/\left(\oplus_{i>n} \mc{L}_i \right).
$$
 We set $\hat{\mc{L}}=\varprojlim_{n} \mc{L}_{\leq n}$. For every field extension
$K\supset \Q_p$ we define 
$\hat{\mc{L}}_K:=\varprojlim_{n} (\mc{L}_{\leq n}\otimes_{\Q_p} K)$.
\end{definition}

We will be only interested in the finite dimensional graded representations 
of $\mc{L}$, which can be identified with the finite dimensional graded
representations of $\hat{\mc{L}}$.

\subsubsection{}
There is natural element $\epsilon\in \hat{\mc{L}}_K$ defined as follows:
\begin{equation}\label{equation-epsilon}
\epsilon:=\sum_{i>0} (a_i\otimes id_K)(\id),  
\end{equation}
with $\id\in K^{\vee}\otimes_{\Q_p}K$ the canonical element. After choosing
a $\Q_p$-basis $v_1,\dots,v_d$ of $K$, we see that
$$
\epsilon= \sum_{i>0}  \sum_{j=1}^{d} a_i(v_j^{\vee})\otimes v_j.
$$

Let $V$ be a finite dimensional graded representation of $\mc{L}$ then 
${\rm exp}(\epsilon)$ is a unipotent automorphism of $V\otimes_{\Q_p}K$.

\begin{proposition}\label{proposition-L-representations-C_eta-equivalence}
The $\otimes$-functor
\begin{align*}
  \Psi:\text{(finite dim.~graded $\mc{L}$-modules)} &\xr{} \mc{C}_{\eta} \\
      V&\mapsto (V,{\rm exp}(\epsilon))
\end{align*}
is an equivalence of categories.
\begin{proof}
  Note first that $\Psi$ is well-defined, because 
$({\rm exp}(\epsilon)-id)(V_n)\subset \oplus_{i>n} V_i$, for all $n$, thus
(\ref{equation-grading-and-eta}) is satisfied.

Let us prove that $\Psi$ is essentially surjective. Let $(\oplus_n V_n,\eta)$
be an object of $\mc{C}_{\eta}$. Since $1-\eta$ is nilpotent we can define
$$
\tilde{\epsilon}=\log(\eta)=\log(1-(1-\eta)).
$$
For every $i>0$ we define a $\Q_p$-linear map $\beta_i:K^{\vee}\xr{} \End(V)_i$
by 
$$
\beta_i(f)=\sum_{n} (id_{V_{n+i}}\otimes f)\circ {\rm proj}_{V_{n+i}\otimes K} \circ 
\tilde{\epsilon} \circ {\rm incl}_{V_n},
$$
where ${\rm incl}_{V_n}:V_n \xr{} V$ and 
${\rm proj}_{V_{n+i}\otimes K}:V\otimes K\xr{} V_{n+i}\otimes K$ is the inclusion
and the projection, respectively. Via (\ref{equation-definition-L}) we
obtain a graded representation $\rho:\mc{L}\xr{} \End(V)$. We need to show that 
${\rm exp}(\rho(\epsilon))=\eta$, or equivalently 
$\rho(\epsilon)=\tilde{\epsilon}$. This is a straight forward computation 
which we leave to the reader.

Next, we need to prove that $\Psi$ is fully faithful. Clearly, $\Psi$ is 
faithful. Let $V,U$ be two graded representations of $\mc{L}$, and suppose
$\tau:V\xr{} U$ is a morphism which respects the grading and commutes 
with $\exp(\epsilon)$. Then $\tau$ commutes with $\epsilon$. Fix a basis 
$v_1,\dots,v_d$ of $K$. For $v\in V$ we get 
\begin{align*}
  \tau\epsilon(v)&=\tau(\sum_{i>0}\sum_j a_i(v^{\vee}_j)(v)\otimes v_j) \\
                 &=\sum_{i>0}\sum_j \tau(a_i(v^{\vee}_j)(v)) \otimes v_j, \\
  \epsilon\tau(v)&=\sum_{i>0}\sum_j a_i(v^{\vee})_j(\tau(v))\otimes v_j.
\end{align*}
Therefore $\tau\circ a_i(v^{\vee}_j) =a_i(v^{\vee}_j) \circ \tau$ for all $i,j$.
Since $\mc{L}$ is generated by the elements $\{a_i(v^{\vee}_j)\}_{i,j}$ we see
that $\tau$ is a morphism of $\mc{L}$-representations.
\end{proof}
\end{proposition}

\begin{corollary}\label{corollary-MT_K-graded-L-modules-equivalent}
  There is an equivalence of $\otimes$-categories 
$$
\Theta:MT_K ^{\phi}\xr{} \text{(finite dim.~graded $\mc{L}$-modules)}
$$
such that 
\begin{itemize}
\item ${\rm forg}\circ \Theta=\tilde{\omega}$, where ${\rm forg}$ forgets
about the $\mc{L}$-action.
\item ${\rm exp}(\epsilon)\mid_{\Theta(M,\phi,F)}=\eta(M,\phi,F)$. 
\end{itemize}
\begin{proof}
Follows immediately from Proposition \ref{proposition-MT_K-Ceta-equivalence}  
and Proposition \ref{proposition-L-representations-C_eta-equivalence}.
\end{proof}
\end{corollary}

\begin{corollary}\label{corollary-MT_K-graded-L-modules-equivalent-2}
  Let $U$ be the pro-algebraic group 
$U=\varprojlim_{n}{\rm exp}(\mc{L}/(\oplus_{i>n} \mc{L}_i))$. Let 
$G_{\tilde{\omega}}$ be the Tannaka group attached to the fibre 
functor $\tilde{\omega}$ (Lemma \ref{lemma-MTK-Tannaka}).

There is an isomorphism
$$
G_{\tilde{\omega}}\cong {\mathbb G}_{m}\ltimes U 
$$ 
such that $\eta\in G_{\tilde{\omega}}(K)$ corresponds to 
${\rm exp}(\epsilon)\in U(K)$.
\begin{proof}
  In view of Corollary \ref{corollary-MT_K-graded-L-modules-equivalent}, the
statement follows from the fact that ${\mathbb G}_{m}\ltimes U$ is the Tannaka 
group of the fibre functor:
$$
{\rm forg}:\text{(graded finite dimensional $\mc{L}$-modules)}\xr{} \text{($\Q_p$-vector spaces)}.
$$
The action of $\mathbb{G}_m$ on $U$ is induced by the action of $\mathbb{G}_m$ on 
$\mc{L}$ given by the grading:
$$
\mathbb{G}_m\times \mc{L}_i \xr{} \mc{L}_i; \quad (t,x)\mapsto t^i\cdot x.
$$
\end{proof}
\end{corollary}

\section{Mixed Tate motives over a number field and logarithmic points}

\subsection{Mixed Tate motives}

\subsubsection{}\label{section-notation-K-S}
Let $E$ be a number field and $S$ a set of finite places. Let $\OO$
be the ring of integers of $E,$ and $|{\rm Spec} (\OO)|$ the maximal spectrum of $\OO.$ 
We denote by 
$$
\OO_S:=\bigcap_{x\in \Spec(\OO)\backslash S}\OO_{x}
$$ 
the ring of  $S$-integers of $E$; the elements of $\OO_S$ are integral outside of $S$. 
We will be mainly interested in two cases for $S$.
In the first case, we have $S=|{\rm Spec} (\OO)|$ and $\OO_S=E$. In the
second case, we have $S=|\Spec(\OO)|\backslash \{x\}$,
for a point $x\in |\Spec(\OO)|$, and $\OO_S=\OO_x$ is the local 
ring at $x$.

\subsubsection{}
Deligne and Goncharov defined in \cite[1.6]{DG} an abelian category
of mixed Tate motives $MT(\OO_S)$. By definition it is the full subcategory
of $MT(E)$ consisting of objects which are unramified outside $S$ in the 
following sense. Let $x\in |\Spec(\OO)|$ be a point lying 
over a prime $p$; then we say that $M\in MT(E)$ is unramified at $x$
if for all primes $\ell\neq p$ the corresponding Galois representation
$M_{\ell}$ is unramified at $x$, i.e.~the inertia subgroup $I_x$ (which is only 
well-defined up to conjugation) 
acts trivially at $M_{\ell}$ \cite[Proposition 1.8]{DG}.

\subsubsection{} 
For extensions of Tate objects we know that: 
\begin{align*}
  \Ext^{1}_{MT(\OO_S)}(\Q(0),\Q(1))&=\OO_{S}^{\times}\otimes \Q,\\
  \Ext^{1}_{MT(\OO_S)}(\Q(0),\Q(n))&=
\begin{cases}
0  &\text{if $n\leq 0$,}\\
\Ext^{1}_{MT(E)}(\Q(0),\Q(n)) &\text{if $n\neq 1$,}
\end{cases}\\
   \Ext^{2}_{MT(\OO_S)}(\Q(0),\Q(n))&=0 \quad \text{for all  $n\in \Z$,}
\end{align*}
(see \cite[Proposition~1.9]{DG}).

\subsubsection{}
Every object of $MT(\OO_S)$ comes equipped with a finite increasing functorial 
weight filtration, 
indexed by even integers. For all $n\in \Z$ the graded pieces 
${\rm gr}^W_{2n}(M)$ are sums of copies of $\Q(-n)$.

In view of \cite[1.1]{DG} the $\otimes$-functor
\begin{align}\label{align-fibre-functor}
\omega: MT(\OO_S)&\xr{} \text{($\Q$-vector spaces)}, \\
\omega(M)&:=\bigoplus_{n\in \Z} \Hom(\Q(n),{\rm gr}^W_{-2n}(M)), \nonumber  
\end{align}
is a fibre functor, therefore $MT(\OO_S)$ is a Tannaka category. We denote
by $G_{S,\omega}$ the group scheme of $\otimes$-automorphisms of $\omega$. 
By \cite[2.1]{DG} we can write $G_{S,\omega}$ as a semi-direct product:
$$
G_{S,\omega}=\mathbb{G}_m\ltimes U_{S,\omega},
$$
where $U_{S,\omega}$ is a unipotent group  and $G_{S,\omega}\xr{} \mathbb{G}_m$ is
induced by the obvious grading of $\omega$. If $S=|{\rm Spec} (\OO)|$ then 
we simply write $G_{\omega}=G_{S,\omega}$.

\subsubsection{Functor to $p$-adic representations}

Let $x\in |\Spec(\OO)|$ be a point lying over a prime $p$. Let $K=E_x$ be the completion of $E$ at the place $x$. Choose algebraic closures
$\bar{E}$, $\bar{K}$, and an embedding $\iota:\bar{E}\xr{} \bar{K}$.

To $M\in MT(E)$ we can attach a Galois representation $M_p$ of 
${\rm Gal}(\bar{E}/E)$ with coefficients in $\Q_p$, which is called the 
$p$-adic realization of $M$. By using $\iota$, we get a continuous homomorphism
$$
{\rm Gal}(\bar{K}/K) \xr{} {\rm Gal}(\bar{E}/E), 
$$ 
and we can restrict $M_p$ in order to obtain a $p$-adic representation 
$M_{\iota,p}$ of ${\rm Gal}(\bar{K}/K)$.

\begin{proposition} \label{proposition-MT(E)-to-MTGK}
  The assignment $M\mapsto M_{\iota,p}$ defines a functor 
$$
(.)_{\iota,p}: MT(E)\xr{} MT_{G_K}.
$$ 
See Definition \ref{definition-mixed-Tate-representation} for $MT_{G_K}$.
\begin{proof}
The $p$-adic realization is functorial. Thus we only need to show that 
$M_{\iota,p}\in MT_{G_K}$, which  follows immediately from the existence of the weight filtration of $M$
and Definition \ref{definition-mixed-Tate-representation}.
\end{proof}
\end{proposition}

The set $\{\iota:\bar{E}\xr{} \bar{E}_{x}\}$ of embeddings over $E$ is
a torsor under the Galois group ${\rm Gal}(\bar{E}/E)$, and for every
$g\in {\rm Gal}(\bar{E}/E)$ there is a natural transformation:
\begin{equation}
\label{equation-iota-to-iotag-1}
\alpha_g:(.)_{\iota,p} \xr{\cong} (.)_{\iota\circ g,p}.  
\end{equation}

\begin{lemma}\label{lemma-comparison-fibre-functors-1}
For the fibre functor $\tau$ (defined in (\ref{definition-tau})) and 
the fibre functor $\omega$
defined in (\ref{align-fibre-functor}) we have a canonical isomorphism
$$
\tau  \circ (.)_{\iota,p} \cong \omega\otimes_{\Q} \Q_p.
$$
For every $g\in {\rm Gal}(\bar{E}/E)$, the diagram 
\begin{equation}\label{diagram-omega-omegatilde-iota-g-1}
\xymatrix
{\tau \circ (.)_{\iota,p}\ar[rd]_{\cong} \ar[rr]^{\text{~(\ref{equation-iota-to-iotag-1})}}
&
&
\tau \circ (.)_{\iota\circ g,p} \ar[ld]^{\cong}
\\
& 
\omega\otimes_{\Q}\Q_p
&
}
\end{equation}
is commutative.
\begin{proof}
Straightforward.
\end{proof}
\end{lemma}

\subsubsection{}
Recall that we have constructed a $K_{st}$-valued $\eta_{st}$ of $\underline{{\rm  Aut}}^{\otimes}\tau$ (Lemma \ref{lemma-eta-st-well-defined}).

\begin{proposition}\label{propostion-eta-st-for-MT(E)-well-defined}
  For every embedding $\iota$, $\eta_x:=\eta_{st}\circ (.)_{\iota,p}$ defines
a $K_{st}$-valued point of $\underline{{\rm  Aut}}^{\otimes}_{MT(E)}\omega$ which 
is independent of the choice of $\iota$.
\begin{proof}
  Since $\tau\circ (.)_{\iota,p}=\omega\otimes_{\Q}\Q_p$ by Lemma \ref{lemma-comparison-fibre-functors-1}, $\eta_{st}\circ (.)_{\iota,p}$ is a $K_{st}$-valued point of $\underline{{\rm  Aut}}^{\otimes}\omega$.

The independence of the choice of $\iota$ follows from the commutative diagram
$$
\xymatrix
{
(\tau \circ (.)_{\iota,p})\otimes_{\Q_p} K_{st} \ar[dd]_{\tau\text{(\ref{equation-iota-to-iotag-1})}\otimes K_{st}}\ar[rrr]^{\eta_{st}} \ar[rd]
&
&
&
(\tau \circ (.)_{\iota,p})\otimes_{\Q_p} K_{st} \ar[dd]^{\tau\text{(\ref{equation-iota-to-iotag-1})}\otimes K_{st}} \ar[ld]
\\
&
\omega\otimes_{\Q}K_{st}
&
\omega\otimes_{\Q}K_{st}
&
\\
(\tau \circ (.)_{\iota\circ g,p})\otimes_{\Q_p} K_{st} \ar[rrr]^{\eta_{st}} \ar[ru]
&
&
&
(\tau \circ (.)_{\iota\circ g,p})\otimes_{\Q_p} K_{st}, \ar[lu]
}
$$
where the triangles are commutative by Lemma \ref{lemma-comparison-fibre-functors-1},
and the square is commutative because $\eta_{st}$ is functorial.
\end{proof}
\end{proposition}

\subsection{Crystalline characterization of unramified motives}

\subsubsection{}
Let $E$ be a number field, and let $M$ be a mixed Tate motive over $E$, i.e.~an 
object in $MT(E)$.     
Let $\nu$ be a finite place of $E,$ $M$ is  unramified at $\nu$  
\cite[Definition 1.4, \textsection 1.7]{DG}  if the coaction \cite[(1.2.2)]{DG} 
$$
e_{M}: \omega(M) \to {\rm Ext}^{1} (\mathbb{Q}(0),\mathbb{Q}(1))\otimes \omega(M)
$$
of ${\rm Ext}^{1} (\mathbb{Q}(0),\mathbb{Q}(1))=E^{\times}\otimes_{\Z}\Q$ on $\omega(M)$ factors through a coaction of $\mathcal{O}_{\nu}^{\times}\otimes_{\Z}\Q.$ 

\subsubsection{} 


Recall from Proposition \ref{proposition-MT(E)-to-MTGK} that
$M_{\iota,p}$ is a mixed Tate Galois representation of 
$G_K={\rm Gal}(\bar{K}/K)$ for the completion $K=E_{\nu}$ at $\nu$. 
In particular, $M_{\iota,p}$ is semistable (Proposition \ref{proposition-mixed-Tate-rep-well-defined}). 
In the following we will simply write $M_p=M_{p,\iota}$. We call $M_p$
\emph{crystalline} if the monodromy operator $N$ of $D_{st}(M_p)$ is trivial, or
equivalently if 
$$
(B_{cris}\otimes_{\Q_p} M_p)^{G_K} \xr{} (B_{st}\otimes_{\Q_p} M_p)^{G_K} 
$$
is an isomorphism.

\begin{thm}\label{thm-unramified=crystalline}
Let $M$ be a mixed Tate motive over $E$ and $\nu$ a finite place of $E.$ Then $M$ is unramified at $\nu$ if and only if $M_{p}$ is crystalline.

\begin{proof}
First note that the statement that $M$  is unramified at $\nu$ is equivalent to the statement that for every subquotient $N$ of $M$
 which is of the form 
 $$
 0 \to \mathbb{Q}(n+1) \to N \to \mathbb{Q}(n) \to 0,
 $$
for some $n,$ the extension class ${\rm Ext}^{1} (\mathbb{Q}(n),\mathbb{Q}(n+1))={\rm Ext} ^{1} (\mathbb{Q}(0),\mathbb{Q}(1))=E ^{\times} \otimes \Q$ lies in $\mathcal{O}_{\nu} ^{\times}\otimes_{\Z}{\mathbb{Q}}$ \cite[\textsection 1.4]{DG}.

Also in the category of $p$-adic representations of a $p$-adic field $K,$  a representation in   ${\rm Ext}^{1} _{G_{K}} (\mathbb{Q}_{p}(0),\mathbb{Q}_{p}(1))$ that is associated to some $q \in K^{\times}\otimes \Q \subseteq \varprojlim _{n} (K^{\times }/ (K^{\times}) ^{p^n})\otimes \mathbb{Q}={\rm Ext}^{1} _{G_{K}} (\mathbb{Q}_{p}(0),\mathbb{Q}_{p}(1))$  is crystalline if and only if $q \in \mathcal{O}^{\times} _{K}\otimes \Q$ \cite[Example 2.3.2]{tsu}. 

First, suppose that $M_{p}$ is crystalline, then every subquotient of $M_{p}$ is crystalline. 
So in order to prove that $M$ is unramified at $\nu$ we may assume that $M=N,$ where 
$N$ is as above with $n=0$ (after Tate twist). 
Therefore we have an extension in ${\rm Ext}^{1} (\mathbb{Q}(0),\mathbb{Q}(1)),$ defined by some $q \in E^{\times}\otimes_{\Z} \Q,$ 
whose $p$-adic realization  is crystalline at $\nu.$  Then the above remark implies that 
the image of $q$ in $E_{\nu} ^{\times}\otimes \Q$ lies in 
$\hat{\mathcal{O}}_{\nu} ^{\times}\otimes \Q$, hence $q \in  \mathcal{O}_{\nu} ^{\times}\otimes \Q$ and $M$ is unramified at $\nu.$ 

Suppose conversely that $M$ is unramified at $\nu.$ 
We have to show that the monodromy operator $N$ on $D_{st}(M_{p})=:D(M)$ vanishes.
Note that $N$ maps the slope  $\lambda$ piece of $D(M)$ to the slope $\lambda -1$ piece. Therefore, if $N$ is nonzero on $D(M)$ then there exists an $n$ 
such that $N$ is nonzero on $D(W_{2n}M_{p} /W_{2n-4}M_{p} )=D((W_{2n}M/W_{2n-4}M)_{p}).$  
Replacing $M$ by $(W_{2n}M/W_{2n-4}M) \otimes \mathbb{Q}(n) $   
we may assume  that $M$ is defined by a class in ${\rm Ext}^{1} (\mathbb{Q}(0) ^{\oplus r},\mathbb{Q}(1) ^{\oplus s})={\rm Ext}^{1} (\mathbb{Q}(0) ,\mathbb{Q}(1))^{\oplus rs}$, $M$ is unramified, and $N$ is nonzero on $D(M)$. 
By passing to a subquotient we may further assume that $r=s=1.$ This gives an extension 
in ${\rm Ext}^{1} (\mathbb{Q}(0),\mathbb{Q}(1))$ which is unramified at $\nu$ (and hence defined by some $q \in \mathcal{O} _{\nu} ^{\times} \otimes \mathbb{Q}$) and whose $p$-adic realization is not  crystalline at $\nu.$ This is a contradiction.  
\end{proof}
\end{thm}


\subsubsection{}
Recall the notation of Section \ref{section-notation-K-S}. Let $x\in |\Spec(\OO)|$
be a point; in the following we will work with $S=|\Spec(\OO)|\backslash \{x\}$,
thus $\OO_S=\OO_x$.

Let $p$ be the prime lying under $x$.
In view of Theorem \ref{thm-unramified=crystalline}, we know that
$MT(\OO_x)$ is the full subcategory of $MT(E)$ 
consisting of motives $M$ such that the $p$-adic realization $M_p$
is crystalline at $x$.

We denote
by $G_{x}$ the group scheme of $\otimes$-automorphisms of the fibre functor (see \eqref{align-fibre-functor})
\begin{equation}\label{equation-omega-local-ring}
\omega:MT(\OO_x) \xr{} \text{($\Q$-vector spaces)}.  
\end{equation}
The group scheme $G_{x}$ is a quotient of $G_{\omega}={\rm Aut}^{\otimes}_{MT(E)}\omega$. 

\begin{lemma}\label{lemma-eta-crys}
The morphism $\Spec(E_{x,st})\xr{\eta_x} G_{\omega}\xr{} G_{x}$ factors through the structure morphism $\Spec(E_{x,st})\xr{} \Spec(E_x)$ and thus defines a point 
$\eta^{ur}_x\in G_x(E_x)$.
\begin{proof}
  The point $\eta_x$ was defined in Proposition
\ref{propostion-eta-st-for-MT(E)-well-defined}. If $M\in MT(\OO_x)$ then 
$D_{st}(M_{\iota,p})$ has vanishing monodromy operator $N$ and 
$\eta_{st}(M_{\iota,p})=\eta_{st,\log}D_{st}(M_{\iota,p})$ takes values in $E_x$
by Definition \ref{definition-eta-st}. 
\end{proof}
\end{lemma}

\subsection{Main theorem}

\subsubsection{}
Let $x\in |\Spec(\OO)|$ and let $E_x$ be the completion of $E$ at $x$.
Bloch and Kato \cite[Definition~3.10]{BK} define an exponential map  
\begin{equation}\label{equation-Bloch-Kato-exp-1}
{\rm exp}:E_x \xr{} \Ext^1(\Q_p(0),\Q_p(n)), \qquad \text{for all $n\geq 1$,}    
\end{equation}
where $\Ext^1$ is computed in the category of $p$-adic representation of 
${\rm Gal}(\bar{E}_x/E_x)$.  Note that, in fact, the image  of the exponential map lies among the crystalline representations  $\Ext^{1}_{crys}(\mathbb{Q}_{p}(0),\mathbb{Q}_{p}(n))$ \cite[Example 3.9]{BK}.  Via $p$-adic Hodge theory,  we obtain a map
\begin{equation}\label{equation-Bloch-Kato-2}
E_x \xr{} \Ext^1_{crys}(\Q_p(0),\Q_p(n)) \cong \Ext^1_{MT^{\phi}_{E_x}}(E_{x}(0),E_{x}(n)),  
\end{equation}
which, by abuse of notation, will also  be called the Bloch-Kato exponential map.

\subsubsection{}
For an extension $M\in \Ext^1_{MT(\OO_x)}(\Q(0),\Q(n))$ with $n\geq 1$, there are
natural maps $v_0:\Q\xr{} \omega(M)$ and $f_n:\omega(M)\xr{} \Q$ defined
as follows. By definition, there are isomorphisms $\alpha:\Q(n)\xr{} {\rm gr}^W_{-2n}M$ and $\beta:{\rm gr}^W_0M\xr{} \Q(0)$; we define 
\begin{align*}
v_0&:\Q = \Hom(\Q(0),\Q(0)) \xr{\beta^{-1}} \omega_0(M) \xr{} \omega(M), \\   
f_n&:\omega(M)\xr{} \omega_n(M) \xr{\alpha^{-1}} \Hom(\Q(n),\Q(n))=\Q. 
\end{align*}
Therefore, we can attach to $M$ a function in $\mathbb{A}^1(G_{x})$, defined by 
$$
M(t):=f_n(t\cdot v_0),
$$
for every point $t:T\xr{} G_x$.

\begin{thm}\label{thm-main}
Let $E$ be a number field and $\OO$ be the ring of integers.  
Let $x\in |\Spec(\OO)|$ be a closed point over a prime $p$. For the Tannaka category $(MT(\OO_x),\omega)$ of 
mixed Tate motives we denote by $G_{x}$ the group scheme of 
$\otimes$-automorphisms of $\omega$.
For all $n\geq 1$, the map 
$$
\Ext^1_{MT(\OO_x)}(\Q(0),\Q(n)) \xr{} E_x, \quad M\mapsto M(\eta^{ur}_x), 
$$
induced by $\eta^{ur}_x\in G_x(E_x)$ (see Lemma \ref{lemma-eta-crys}),
is the composition of the p-adic realization 
$$
\Ext^1_{MT(\OO_x)}(\Q(0),\Q(n)) \xr{} \Ext^1 _{crys}(\Q_p(0),\Q_p(n))
$$ 
and 
the inverse of the Bloch-Kato exponential map (\ref{equation-Bloch-Kato-exp-1}).

\begin{proof}

Let us prove that evaluation at the point $\eta^{ur}_x$ has the desired compatibility
with the Bloch-Kato exponential map (see (\ref{equation-Bloch-Kato-2}))
$$
{\rm exp}:E_x\xr{} \Ext^1 _{crys}(\Q_p(0),\Q_p(n)) \xr{\cong} \Ext^1_{MT_{E_x}}(E_{x}(0),E_{x}(n)). 
$$
For this we need to recall the construction of 
the exponential map. For the rest of the proof let $K:=E_{x}.$   First there is an exact sequence \cite[Proposition 1.17]{BK}:
\begin{equation}\label{exact seq bk}
0 \to \mathbb{Q}_{p} \to B_{crys} ^{\varphi=1} \oplus B_{dR} ^{+} \to B_{dR}\to 0,
\end{equation}
where the first map sends $x$ to $(x,x)$ and the second one sends $(x,y)$ to $x-y.$ 

For $n \geq 1,$ the Bloch-Kato construction gives a map 
$$
K=(\mathbb{Q}_{p}(n) \otimes B_{dR})^{G_{K}} \to {\rm Ext} ^{1} _{crys} (\mathbb{Q}_{p}(0), \mathbb{Q}_{p}(n)).
$$ 
This map is obtained as follows. First tensor the above exact sequence with $\mathbb{Q}_{p}(n)$:
$$
0\to \mathbb{Q}_{p}(n) \to (\mathbb{Q}_{p}(n) \otimes B_{crys} ^{\varphi=1}) \oplus (\mathbb{Q}_{p}(n) \otimes B_{dR} ^{+}) \to \mathbb{Q}_{p}(n) \otimes B_{dR}\to 0.
$$

Then an element $a$ in $K(n)=(\mathbb{Q}_{p}(n)\otimes B_{dR})^{G_{K}}$ gives a map $\mathbb{Q}_{p} \to \mathbb{Q}_{p}(n)\otimes B_{dR},$ pulling back the above exact sequence via this map gives the extension we were looking for. 

More explicitly, for $a \in K$ the extension constructed above is:
$$
0 \to V_{n} \to V \to V_{0} \to 0,
$$
where $V_{0}=\mathbb{Q}_{p}\cdot  t^{n} \otimes a t^{-n},$ $V_{n}=\mathbb{Q}_{p} \cdot t^{n},$ and $V$ is a 2-dimensional representation of $G_{K}$ with basis which can be described as follows. 
By the exact sequence  (\ref{exact seq bk}), there exists $x \in B_{crys} ^{\varphi=1}$ and $y \in B_{dR} ^{+}$ such that $at^{-n}=x-y.$  
Then $V$ has basis $\{ (t^{n} \otimes x,t^{n} \otimes y), (t^{n} \otimes 1, t^{n} \otimes 1)  \}. $    For $\sigma \in G_{K},$ 
$$
\sigma (t^{n} \otimes x,t^{n}\otimes y)=(t^{n} \otimes x,t^{n}\otimes y)+\gamma(\sigma) (t^{n} \otimes 1, t^{n} \otimes 1),
$$
for some $\gamma(\sigma) \in \mathbb{Q}_{p}.$ Therefore 
$$
\chi_{cyc} (\sigma) ^{n} \sigma (x)=x+\gamma(\sigma)
$$
and 
$$
\chi_{cyc} (\sigma) ^{n} \sigma (y)=y+\gamma(\sigma)
$$

Let us now try to find what this extension corresponds to after we apply the functor $(\cdot \otimes B_{crys})^{G_{K}}.$ First note that $(V\otimes B_{crys})^{G_{K}}$ has basis 
$$
e_{n}:= (t^{n}\otimes 1,t^{n} \otimes 1)\otimes t^{-n}
$$
and 
$$
e_{0}:=(t^{n} \otimes x,t^{n} \otimes y)\otimes 1 -(t^{n} \otimes 1, t^{n} \otimes 1) \otimes x.
$$
That $e_{n}$ is invariant under the Galois action is clear.  
In order to see that $e_{0}$  is $G_{K}$ invariant let $\sigma \in G_{K}.$  
Then 
$$
\sigma (e_{0})=(t^{n}\otimes (x+\gamma(\sigma)), t^{n} \otimes (y + \gamma(\sigma)) )\otimes 1-(t^{n}\otimes 1, t^{n} \otimes 1)\otimes (x+ \gamma(\sigma))=e_{0}.
$$ 
Now note that $\varphi(e_{n}) =p^{-n}e_{n}$ and $\varphi(e_{0})=e_{0}.$ Furthermore $e_{n}$ is the image of $1 \in K(n)$  and $e_{0}$ maps to $1 \in K(0)$ in the exact sequence (note that $\mathbb{Q}_{p}(0)$ is identified with $V_{0}$ via the map that sends 1 to $t^{n} \otimes a t^{-n}$): 
$$
0 \to K(n) \to (V\otimes B_{crys})^{G_{K}} \to K(0) \to 0.
$$
Therefore in order to compare Bloch-Kato's construction we need only compute the filtration on $(V\otimes B_{crys})^{G_{K}} \otimes _{K_{0}}K.$ So we need to compute the 0-th piece of the filtration on $(V \otimes B_{dR})^{G_{K}}.$ 
  
We claim that $ae_{n} +e_{0} \in Fil ^{0} (V \otimes B_{dR})^{G_{K}}.$ This follows immediately from 
$$
ae_{n}+e_{0} = (t^{n} \otimes x,t^{n} \otimes y)\otimes 1-(t^{n} \otimes 1,t^{n} \otimes 1)\otimes y,
$$
and the fact that $y \in B_{dR}^{+}.$ Now Proposition \ref{proposition-MTK-extensions} 
implies the claim.
\end{proof}
\end{thm}

\subsection{Archimedian places}

In this section we recall the story for archimedian places; our reference 
is \cite{D} and \cite{BD}.

\subsubsection{}
Let $E$ be a number field and $\sigma:E\xr{} \C$ an embedding. To 
$M\in MT(E)$ we can attach a \emph{real} mixed Tate Hodge structure $M_{\sigma}$.
Recall that a real mixed Tate Hodge structure $(H,W,F)$ consists of an $\R$-vector 
space $H$, an increasing filtration $W$ of $H$, and a decreasing filtration $F$ of $H\otimes_{\R}\C$ such that  
$$
{\rm Gr}^p_{F}{\rm Gr}^q_{\bar{F}}{\rm Gr}^W_{n} (H\otimes \C)=
\begin{cases}
0 &\text{if $n$ is odd,} \\
0 &\text{if $n$ is even and $(p,q)\neq (\frac{n}{2},\frac{n}{2})$.}
\end{cases} 
$$

Induced by $F,\bar{F}$, we obtain maps 
\begin{align*}
  a_{F}: H\otimes \C &= \bigoplus_{i\in \Z} F^{-i}\cap W_{-2i} \xr{} \bigoplus_{i\in \Z} {\rm Gr}^W_{-2i}H\otimes \C, \\
 a_{\bar{F}}: H\otimes \C &= \bigoplus_{i\in \Z} \bar{F}^{-i}\cap W_{-2i} \xr{} \bigoplus_{i\in \Z} {\rm Gr}^W_{-2i}H\otimes \C,
\end{align*}
where $F^{i}\cap W_{2i}\xr{} {\rm Gr}^W_{2i}H\otimes \C$ is the natural map (and similarly 
for $\bar{F}$). For the automorphism $d=a_{\bar{F}}a_F^{-1}$ of 
$\bigoplus_{i\in \Z} {\rm Gr}^W_{-2i}H\otimes \C$ we know that
$$
(d-1)({\rm Gr}^W_{-2i}H\otimes \C)\subset \bigoplus_{j>i} {\rm Gr}^W_{-2j}H\otimes \C,
$$ by  \cite[p.510]{D}, and $\bar{d}=d^{-1}$ \cite[p.513]{D}.

\subsubsection{}
Let $\mc{C}$ be the category of pairs $(\oplus_i H_i,d)$ where $\oplus_i H_i$
is a graded $\R$-vector space and 
$d:\oplus_i H_i\otimes \C\xr{} \oplus_i H_i\otimes \C$ is an automorphism 
satisfying the conditions $\bar{d}=d^{-1}$ and $(d-1)(H_i)\subset \oplus_{j> i} H_{j}$, for all 
$i$.

\begin{proposition}\cite[p.514]{D}
  The functor 
  \begin{align*}
    \text{(Real mixed Tate Hodge structures)}\xr{} \mc{C} \\
    (H,W,F)\mapsto (\oplus_{i\in \Z} {\rm Gr}^W_{-2i}H,d)
  \end{align*}
is an equivalence of categories.
\end{proposition}

The maps $d$ define a $\C$-valued $\otimes$-automorphism for the fibre functor
\begin{align*}
&\tilde{\omega}:\text{(Real mixed Tate Hodge structures)} \xr{} \text{($\R$-vector spaces)},  \\
&\tilde{\omega}(H,W,F)=\bigoplus_{i\in \Z}{\rm Gr}^W_{-2i}H.
\end{align*}

\subsubsection{}
Recall the definition of $\omega$ in (\ref{align-fibre-functor}). 
For the functor
$$
\mc{R}_{\sigma}: MT(E)\xr{} \text{(Real mixed Tate Hodge structures)}, \quad 
M\mapsto M_{\sigma},
$$  
we have an isomorphism
\begin{equation}\label{equation-omega-R}
\omega\otimes_{\Q}\R \cong \tilde{\omega}\circ \mc{R},  
\end{equation}
depending on the choice $(2\pi i)^n$ as a generator for the real vector space 
underlying $\R(n)$, in other
words we have to choose a square root of $-1$ in $\C$. In order to avoid 
this choice one can define 
$$
\tilde{\omega}(H,W,F)=\bigoplus_{n\in \Z} i^n \cdot {\rm Gr}^W_{-2n}H,
$$
as in \cite[p.111]{BD}, but we won't do that. 

Via Equation (\ref{equation-omega-R}), $d$ defines a $\C$-valued point 
of $G_{\omega}$, the $\otimes$-automorphisms of the fibre functor $\omega$.
We define $\epsilon=\log(d)$, $\epsilon$ defines a $\C$-valued point of 
${\rm Lie}(G_{\omega})$. 

The dictionary for the notation of \cite[p.111]{BD} is
$$
d=b^{-1}, \quad \epsilon=-2\cdot N,
$$
and $N$ is purely imaginary.

\subsubsection{}
For $z\in E\backslash \{0,1\}$ there is a polylogarithm motive 
$\{z\}\in MT(E)$ (strictly speaking it is a pro-object). The motive 
$\{z\}$ is defined as a subquotient of the motivic paths from the 
tangent vector $t_0=z$, in the tangent space at $0$, to $z$ (see \cite[Theorem~4.4]{DG}). 
The $\Q$-Hodge realization of $\{z\}$ is described in \cite[p.98]{BD}
and uniquely determines $\{z\}$.

For every $k\in \Z_{\geq 0}$ we have natural isomorphisms
$$
\alpha_k:{\rm gr}^W_{-2k}\{z\}\xr{\cong} \Q(k);
$$ 
we define $v_0\in \omega(\{z\})$  and $f_k\in \omega(\{z\})^{\vee}$ by 
\begin{align*}
v_0:\Q = \Hom(\Q(0),\Q(0)) &\xr{\alpha_0^{-1}} \omega_0(\{z\}) \\   
f_k:\omega(\{z\})\xr{} \omega_k(\{z\}) &\xr{\alpha_k} \Hom(\Q(k),\Q(k))=\Q. 
\end{align*}
We denote by $(v_0,\{z\},f_k)\in \mathbb{A}^1({\rm Lie}\; G_{\omega})$ the function 
$$
X\mapsto f_k(X\cdot v_0).
$$

By \cite[Proposition~2.7]{BD} we have
$$
(v_0,\{z\},f_k)(\epsilon)=
\begin{cases}
2i\sum_{\ell} b_\ell \frac{{\rm log}(z\bar{z})}{\ell!} {\rm Im}({\rm Li}_{k-\ell}(z)) &\text{if $k$ is even,}\\
2i\sum_{\ell} b_\ell \frac{{\rm log}(z\bar{z})}{\ell!} {\rm Re}({\rm Li}_{k-\ell}(z))
&\text{if $k$ is odd.}
\end{cases}
$$
Here, $\{b_{\ell}\}$ are the Bernoulli numbers and ${\rm Li}$ is the polylogarithm.
For $k$ even, the result does not depend on the choice of the square root
of $-1$; for $k$ odd it is independent up to a sign.









\end{document}